\documentclass[12pt,a4paper,leqno,dvips,twoside]{article}

\usepackage{graphicx}
\usepackage{amsfonts}
\usepackage{amssymb}
\usepackage{amsmath}
\usepackage{amsthm}

 \numberwithin{equation}{section}
\newtheorem{theo}{Theorem}[section]

\newtheorem{defin}{Definition}[section]
\newtheorem{lema}{Lemma}[section]

\newtheorem{cor}{Corollary}[section]

\newtheorem{rema}{Remark}[section]

\newtheorem{example}{Example}

\begin{document}

\setcounter{page}{1}
\title{ A note about the relation between  fixed point theory on cone metric spaces and
fixed point theory on  metric spaces}
\author{ Ion Marian Olaru}
\maketitle

\begin{abstract}
  Let $Y$ be a locally convex Hausdorff  space, $K\subset E$ a cone and $\leq_K$
   the partial order  defined by $K$. Let  $(X,p)$ be a $TVS-$ cone metric space,
    $\varphi:K\rightarrow K$  a  vectorial comparison function and $f:X\rightarrow X$ such that
  \[p(f(x),f(y))\leq_K\varphi(p(x,y)),\] for all $x,y\in X$. We
  shall show   that there exists a  scalar comparison function $\psi:\mathbb{R}_+\rightarrow
  \mathbb{R}_+$  and a metric $d_p$(in usual sense) on $X$ such that
  \[d_p(f(x),f(y))\leq\psi(d_p(x,y)),\] for all $x,y\in X$. Our
  results extend the results of Du (2010) [Wei-Shih Du, A note on cone metric fixed point theory and its
equivalence, Nonlinear Anal. 72 (2010), 2259-2261].

\end{abstract}

\begin{center} 2010 Mathematical Subject Classification: 47H10,
54H25

Keywords:  K-metric spaces, cone metric space,  TVS- cone metric
spaces, comparison function\end{center}
\section{Introduction and preliminaries}

Fixed point theory in  K-metric and K-normated spaces was
developed by A.I. Perov and his  consortiums (\cite{perov1},
\cite{perov2}, \cite{perov3}). The  main idea consists to use an
ordered Banach space instead of the set of real numbers, as the
codomain for a metric. For more details on fixed point theory in
K-metric and K-normed spaces, we refer the reader to
\cite{zabreiko}. Without mentioning these previous works, Huang
and Zhang \cite{guang} reintroduced such  spaces under the name of
cone metric spaces but went further, defining convergent and
Cauchy sequences in the terms of interior points of the underlying
cone. They also proved some fixed point theorems in such spaces in
the same work. After that, fixed point results in cone metric
spaces have been studied by many other  authors. References
\cite{abas},\cite{radenovici}, \cite{racocevici},
\cite{radenovic}, \cite{raja}, \cite{rezapour} are some works in
this line of research. However, very recently Wei-Shih Du in
\cite{du} used the scalarization function and investigated the
equivalence of vectorial versions of fixed point theorems in cone
metric spaces and scalar versions of fixed point theorems in
metric spaces. He showed that many of the fixed point results in
ordered K-metric spaces  for maps satisfying contractive
conditions of a linear type in K-metric spaces can be considered
as the corollaries of corresponding theorems in metric spaces.

 Let $E$ be a topological vector space (for short $t.v.s$) with its zero vector  $\theta_E$.
\begin {defin}( \cite{du}, \cite{guang}) A subset K of E is called a cone if:
\begin{itemize}
\item[(i)] K is closed, nonempty and $K\neq\{\theta_E\}$;

\item[(ii)]  $a,b\in\mathbb{R}$, $a,b\geq 0$ and $x,y\in K$ imply
$ax+by\in K$;

\item[(iii)]$K\cap -K=\{\theta_E\}.$
\end{itemize}
\end{defin}

For a given cone $K\subset E$, we can define a partial ordering
$\leq_K$ with respect to $K$ by
\begin{equation}
x\leq_K y\  if\  and\  only\  if \ y-x\in K.
\end{equation}
We shall write $x<_Ky$ to indicate that $x\leq_K y$ but $x\neq y$,
while $x\ll y$ will stand for $y-x \in int K$ (interior of $K$).

In the following, unless otherwise specified, we always suppose
that $Y$ is a locally convex Hausdorff  with its zero vector
$\theta$, $K$ a  cone in $Y$ with $int K \neq \emptyset$ , $e\in
int K$ and $\leq_K$ a partial ordering with respect to $K$.

 \begin{defin}\label{def1}( \cite{du})
 Let $X$ be a nonempty set. Suppose that a mapping $d:X\times X\rightarrow
 Y$ satisfies:
 \begin{itemize}
 \item[(i)] $\theta \leq_K d(x,y)$ for all $x,y\in X$ and $d(x,y)=\theta $ if
 and only if $x=y$;

 \item[(ii)] $d(x,y)=d(y,x)$, for all $x,y\in X$ ;

 \item[(iii)] $d(x,y)\leq_K d(x,z)+d(z,y)$  for all $x,y,z\in X$.
 \end{itemize}
  Then $d$ is called a TVS-cone metric  on $X$ and $(X,d)$ is called a
 TVS-cone metric   space.
  \end{defin}

 The nonlinear scalarization function $\xi_e: Y\rightarrow \mathbb{R}$
is defined as follows
\[\xi_e(y)=\inf\{r\in\mathbb{R}\mid  y\in r\cdot e-K\}.\]

\begin{lema}\label{l1}( \cite{chen}) For each $r\in \mathbb{R}$ and $y\in Y$,
the following statements are satisfied:
\begin{itemize}
\item[(i)]$\xi_e(y)\leq r$ if and only if $y\in r\cdot e-K$;

\item[(ii)]$\xi_e(y)> r$ if and only if  $y\notin r\cdot e-K$;

\item[(iii)]$\xi_e(y)\geq r$ if and only if  $y\notin r\cdot e-int
K$;

\item[(iv)]$\xi_e(y)< r$ if and only if $y\in r\cdot e- intK$;

\item[(vi)] $\xi_e(\cdot)$ is positively homogeneous and
continuous on Y;

\item[(vii)] if $y_1\in y_2+K$ then $\xi_e(y_2)\leq \xi_e (y_1)$;

\item[(viii)]$\xi_e(y_1+y_2)\leq \xi_e(y_1)+\xi_e(y_2)$, for all
$y_1, y_2\in Y$.
\end{itemize}
\end{lema}
\begin{theo}( \cite{du}) Let $(X,p)$ be a $TVS-$cone metric space.
Then \[d_p:X\times X\rightarrow [0,\infty)\] defined by $d_p=\xi_e\circ d$ is a metric.\end{theo}
\section{Main results}
\begin{defin}\label{def2} Let $K\subset Y$ be a cone.
A function $\varphi:K\rightarrow K$ is called a vectorial
comparison function if
\begin{itemize}

\item[(i)]$k_1\leq_P k_2$ implies
$\varphi(k_1)\leq_P\varphi(k_2)$;

\item[(ii)] $\varphi(0)=0$ and $0<_P \varphi(k)<_P k$ for  $k\in
K-\{0\}$;

\item[(iii)] $k\in int K$ implies $k-\varphi(k)\in  int K$;

\item[(iv)]  if $t_0\geq 0$ then $\lim\limits_{t\rightarrow
t_0^+}\varphi(t\cdot e)=\varphi(t_0\cdot e)$.
\end{itemize}\end{defin}
\begin{example}
\begin{itemize}
\item[(i)] if K is an arbitrary cone in a Banach space E and
$\lambda\in (0,1)$, then $\varphi:K\rightarrow K$, defined by
$\varphi(k)=\lambda k$ is a  vectorial comparison function;

 \item[(ii)]Let $E=\mathbb{R}^2$, $K=\{(x,y)\mid x,y\geq 0\}$ and
 let  $\varphi_1,\varphi_2:[0,\infty)\rightarrow  [0,\infty)$ be such that

 \begin{itemize}

 \item[(a)] $\varphi_1,\varphi_2$ are increasing functions;
 \item[(b)]  if $t>0$ then $\varphi_i(t)< t$ for $i=\overline{1,2}$;
 \item[(c)]  $\varphi_1,\varphi_2$ are right continuous.
 \end{itemize}
 Then $\varphi:K\rightarrow K$, defined by  $\varphi(x,y)=(\varphi_1(x),\varphi_2(y))$ is a
vectorial comparison function;

\end{itemize}\end{example}

\begin{defin} ( \cite{rus}) A function  $\varphi: \mathbb{R}_+\rightarrow \mathbb{R}_+$  is called a  scalar comparison function
if

\begin{itemize}

\item[(i)] $t_1\leq t_2$ implies $\varphi(t_1)\leq\varphi(t_2)$;

\item[(ii)]$\varphi^n(t)\stackrel{n\rightarrow\infty}{\rightarrow}
0$ for all $t>0$
\end{itemize}\end{defin}
The following lemma will be useful in the sequel
\begin{lema}\label{l3}( \cite{rus}) If $\varphi: \mathbb{R}_+\rightarrow \mathbb{R}_+$ is increasing and right upper
semicontinuous then the following assertions are equivalent:
\begin{itemize}

\item[(a)] $\varphi^n(t)\stackrel{n\rightarrow\infty}{\rightarrow}
0$ for all $t>0$;

\item[(b)] $\varphi(t)<t$ for all $t>0$.
\end{itemize}
\end{lema}

\begin{lema}\label{l2} We consider $M:\mathbb{R}\rightarrow Y, M(r)=r\cdot e $. Then we have

\begin{itemize}
\item[(i)] $M(0)=\theta$;

\item[(ii)] if $r_1\leq r_2$ then $M(r_1)\leq_K M(r_2)$;

\item[(iii)] $y\leq_K M\circ \xi_e(y)$ for all $y\in Y$;

\item[(iv)] $\xi\circ M(r)\leq r$ for all $r\in\mathbb{R}$;

\item[(v)] if $y_1\ll y_2$ then $\xi_e(y_1)<\xi_e(y_2)$.
\end{itemize}
\end{lema}

\noindent{\bf Proof:}

$(i)$ It is obvious;

$(ii)$ Let be $r_1\leq r_2$. Then $(r_2-r_1)\cdot e\in K$. Thus
$M(r_1)\leq_K M(r_2)$;

$(iii)$  Since $\xi_e(y)=\inf\{r\in\mathbb{R}\mid  y\leq_K r\cdot
e\}$ it follows that $y\leq_K \xi_e(y)\cdot e=M\circ\xi_e(y)$ for
all $y\in Y$;

$(iv)$ Let be $r\in\mathbb{R}$. Since $\{r^\prime\in\mathbb{R}\mid
r\cdot e\leq_K r^\prime\cdot e \}\supseteq
\{r^\prime\in\mathbb{R}\mid r\leq r^\prime\}$ we get
\[\xi_e(M(r))=\xi_e(r\cdot e)=\inf\{r^\prime\in\mathbb{R}\mid r\cdot e\leq_K r^\prime\cdot e \}\leq
\inf\{r^\prime\in\mathbb{R}\mid r\leq r^\prime\}=r.\]

$(v)$ Let be  $y_1\ll y_2$.  We remark that  $y_1\ll
y_2\leq_K\xi_e(y_2)\cdot e$. Then, via Remark 1.3 of Radenovi\'{c}
and Kadelburg  \cite{radenovic}, it follows that $y_1\ll
\xi_e(y_2)\cdot e$. Hence $y_1\in \xi_e(y_2)\cdot e- int K$. By
using Lemma \ref{l1} (iv) we get $\xi_e(y_1)<\xi_e(y_2)$.

\begin{theo}\label{t1} Let $(X,p)$ be a TVS-cone metric and $\varphi:K\rightarrow K$
 be a  vectorial comparison function such that

\[p(f(x),f(y))\leq_K \varphi(p(x,y)),\] for all $x,y\in X$. Then there exists a  scalar comparison function
$\psi:\mathbb{R}_+\rightarrow \mathbb{R}_+$ such that

\[d_p(f(x),f(y))\leq \psi(d_p(x,y)),\] for  all $x,y\in X$.
\end{theo}
\noindent{\bf Proof:} Let be $t\in\mathbb{R}_+$. Then
$\theta\leq_K M(t)$. It follows that $M(t)\in K$ for all
$t\in\mathbb{R}_+$.

We define
\[\psi:\mathbb{R}_+\rightarrow \mathbb{R}_+,\]
\[\psi(t)=\xi_e\circ\varphi\circ M(t)\]
First, we note that for all $t\in\mathbb{R}_+$ we have

\[0\leq  \xi_e\circ\varphi\circ M(t)\leq \xi_e\circ
M(t)\leq t.\]

 Now, we remark that for each $x,y\in X$ we have
\[d_p(f(x),f(y))\leq \xi_e\circ \varphi(p(x,y))\leq \xi_e\circ\varphi(M(\xi_e(p(x,y)))) =\psi(d_p(x,y)).\]

We claim that $\psi$ is a scalar  comparison function. Since
$\xi_e$, $\varphi$ and $M$ are increasing functions, it follows
that $\psi $ is increasing function. In order to prove that
$\psi^n(t)\stackrel{n\rightarrow\infty}{\rightarrow} 0$ for all
$t>0$, we shall use  Lemma \ref{l3}.  Next we show that
$\psi(t)<t$ for all $t>0$.

Let be $t_0>0$. Then $t_0\cdot e\in int K$. Therefore
$\varphi(t_0\cdot e)\ll t_0\cdot e$. It follows that
\[\psi(t_0)=\xi_e\circ\varphi(t_0\cdot e)<\xi_e\circ M(t_0)\leq t_0.\]

Since $\lim\limits_{t\rightarrow
t_0^+}\psi(t)=\lim\limits_{t\rightarrow
t_0^+}\xi_e\circ\varphi(t\cdot e)=\xi_e(\lim\limits_{t\rightarrow
t_0^+}\varphi(t\cdot e))=\xi_e\circ\varphi(t_0\cdot e)=\psi(t_0)$
 it follows that $\psi$ is right upper semicontinuous. Hence
 $\psi^n(t)\stackrel{n\rightarrow\infty}{\rightarrow}0$.
\begin{cor}\label{c1} Let $(X,p)$ be a  complete TVS cone metric space and $\varphi:K\rightarrow K$
a vectorial comparison function such that
\[p(fx,fy)\leq_K \varphi(p(x,y)),\]
for all $x,y\in X$. Then, f has a unique fixed point $x_0$.
\end{cor}
\noindent{\bf Proof:} We apply  Theorem \ref{t1} and Theorem 1 pp
459 of Boyd and Wong (\cite{boyd}).
\begin{rema}\label{r1} For   $\varphi(k)=\lambda\cdot k$, $\lambda\in[0,1)$  we obtain, via Lemma \ref{l2} $(iv)$ and Corollary \ref{c1},
the results of W.S. Du \cite{du}.\end{rema}

\begin{rema}\label{r2} Let $(X,p)$ a cone metric space.  For   $\varphi(k)=\lambda\cdot k$, $\lambda\in[0,1)$  we obtain,
via Remark \ref{r1}, the results of L.G. Huang and  Zhang Xian
\cite{guang}.\end{rema}

 Let $(X,d)$ be a {\it TVS} cone-metric
space and let $\varphi:K\rightarrow K$ be a vectorial comparison
function. For a pair $(f,g)$ of self-mappings on $X$ consider the
following conditions:
\begin{itemize}
\item[(C)] for arbitrary $x,y\in X$ there exists $u\in
\{d(gx,gy),d(gx,fx), d(gy,fy)\}$ such that $d(fx,fy)\leq_P
\varphi(u)$.

\item[$(C_1)$] for arbitrary $x,y\in X$ there exists $w\in
\{d_p(gx,gy),d_p(gx,fx), d_p(gy,fy)\}$ such that $d_p(fx,fy)\leq
\psi(u)$.
\end{itemize}
\begin{rema}\label{r3} The condition (C) imply the condition
$(C_1)$.\end{rema}

Indeed since the condition $(C)$ hold, it follows  that  at least
one of the following three cases holds:

\begin{itemize}
\item[\it Case 1:] $u=d(gx,gy)$. Then

 \[\xi_e(p(fx,fy))\leq
\xi_e\circ\varphi(u)\leq \xi_e\circ \varphi\circ
M(\xi_e(u))=\psi(d_p(gx,gy))
\]

\item[\it Case 2:] $u=d(gx,fx)$. Then

 \[\xi_e(p(fx,fy))\leq
\xi_e\circ\varphi(u)\leq \xi_e\circ \varphi\circ
M(\xi_e(u))=\psi(d_p(gx,fx))
\]

\item[\it Case 3:] $u=d(gy,fy)$. Then

 \[\xi_e(p(fx,fy))\leq
\xi_e\circ\varphi(u)\leq \xi_e\circ \varphi\circ
M(\xi_e(u))=\psi(d_p(gy,fy))
\]
\end{itemize}

\noindent Departament of Mathematics,\\
Faculty of Sciences,\\
University "Lucian Blaga" of Sibiu,\\
Dr. Ion Ratiu 5-7, Sibiu, 550012, Romania \\
E-mail: marian.olaru@ulbsibiu.ro

\end{document}